\documentclass[11pt]{amsart}
\usepackage{amsthm}
\usepackage{amsfonts}
\usepackage{amsmath}
\usepackage{amssymb}
\usepackage[all]{xy}

\newtheorem{thm}{Theorem}[section]
\newtheorem{lemma}[thm]{Lemma}
\newtheorem{prop}[thm]{Proposition}

\theoremstyle{definition}
\newtheorem{defn}[thm]{Definition}

\def\ve{\varepsilon}

\def\R{{\mathbb R}}

\def\C{{\mathbb C}}

\def\Z{{\mathbb Z}}

\def\so{_{\scriptscriptstyle O}}
\def\su{_{\scriptscriptstyle U}}
\def\st{_{\scriptscriptstyle T}}
\def\pr{^{\scriptstyle \R}}
\def\po{^{\scriptscriptstyle O}}
\def\pu{^{\scriptscriptstyle U}}
\def\pt{^{\scriptscriptstyle T}}
\def\sr{_{\scriptscriptstyle \R}}
\def\sc{_{\scriptscriptstyle \C}}
\def\crt{^{\scriptscriptstyle {\it CRT}}}
\def\scrt{_{\scriptscriptstyle {\it CRT}}}
\def\ct{{\it CRT}}

\def\im{\text {image} \,}

\def\hom{\text {Hom}}

\everymath{\displaystyle}
\CompileMatrices

\title{The Range of United $K$-theory} 

\date{\today}
\author{Jeffrey L. Boersema}
\address{Department of Mathematics \\
Seattle University \\
Seattle, WA  98122}
\email{boersema@seattleu.edu}

\begin{document}

\begin{abstract}
We prove that the united $K$-theory functor is a surjective functor from the category 
of real simple purely infinite C*-algebras to the category of countable acyclic \ct-modules.
\end{abstract}

\maketitle


\section{Introduction}

In this paper, we further investigate the united $K$-theory 
functor for real C*-algebras, developed in \cite{Boer2}.  
We will show that united $K$-theory is a surjective function from the category of real simple purely infinite C*-algebras to the category of  
countable acyclic \ct-modules, providing one part of a classification-type theorem 
for real C*-algebra along the lines of that of Kirchberg \cite{Kirchberg} 
and of Phillips \cite{Phillips}.  Of course the other part, injectivity, is the more 
difficult aspect to such a classification program, but the universal 
coefficent theorem of \cite{Boer3} described below and the surjectivity 
result of the present paper give us confidence that such a classification 
theorem for real C*-algebras is possible.

The immediate purpose of the development of united $K$-theory was to
state and prove a K\"unneth-type formula for the tensor product 
to two real C*-algebra (in \cite{Boer2}).  Later, in \cite{Boer3}, we proved a universal 
coefficient theorem (UCT) for united $K$-theory (more on this later).  
Recall that for a real C*-algebra $A$, 
united $K$-theory $K\crt(A)$ consists of three graded modules and the 
collection of natural transformations between them.  The three objects are 
\begin{enumerate}
\item[(1)] real $K$-theory $KO_*(A)$ --- defined to be the $K$-theory of 
the real C*-algebra $A$ as discussed for example in \cite{Schroder}.  
\item[(2)] complex $K$-theory $KU_*(A)$ --- defined to be the $K$-theory of the 
complexification $A\sc = \C \otimes A$. 
\item[(3)] self-conjugate $K$-theory $KT_*(A)$ --- defined to be the $K$-theory 
of $T \otimes A = \{ f \colon [0,1] \rightarrow \C \otimes A \mid f(0) = 
\overline{f(1)} \}$.  
\end{enumerate}
These objects are taken not just as graded groups, but as graded modules over the graded unital rings $KO_*(\R)$, $KU_*(\R)$, and $KT_*(\R)$ respectively; which are displayed in degrees $0$ through $8$ by  
\begin{center} \begin{tabular}{c@{$~~=~~$}ccccccccc}
$K_*(\R)$ & $\Z$ & $\Z_2 $ & $\Z_2$ & $0$ & $\Z$ & $0$ & $0$ & $0$ & $\Z$  \\
$K_*(\C)$ & $\Z$ &$0$ & $\Z$ &$0$ & $\Z$ & $0$ &$\Z$& 0 & $\Z$ \\
$K_*(T)$ & $\Z$ & $\Z_2$ & $0$ & $\Z$ & $\Z$ & $\Z_2$ & $0$ & $\Z$ & $\Z$ \\
\end{tabular}   \end{center} 
(see page 23 in \cite{Schroder} and Tables~1, 2, and 3 in \cite{Boer2}).  The 
generators are the elements $1\so \in K_0(\R)$, $\eta\so \in K_1(\R)$, $\eta\so^2 \in K_2(\R)$, $\xi \in K_4(\R)$, 
and the invertible element $\beta\so \in K_8(\R)$.
The ring $K_*(\C)$ is the free unital polynomial ring generated by the invertible Bott element $\beta\su \in K_2(\C)$.
The ring $K_*(T)$ has generators
$1\st$in degree 0, $\eta\st$ in degree 1, $\omega$ in degree 3, and the invertible element $\beta\st$ 
is degree 4.  Thus, $KO_*(A)$ has period 8,  $KU_*(A)$ has period 2, and $KT_*(A)$ has period 4.
United $K$-theory also incorporates the natural transformations  
\begin{align*}   \label{operations}
c_n &\colon KO_n(A) \longrightarrow KU_n(A)  &
r_n &\colon KU_n(A) \longrightarrow KO_n(A) \\
\varepsilon_n &\colon KO_n(A) \longrightarrow KT_n(A) &
\zeta_n &\colon KT_n(A) \longrightarrow KU_n(A) \\
(\psi\su)_n &\colon KU_n(A) \longrightarrow KU_n(A)  &
(\psi\st)_n &\colon KT_n(A) \longrightarrow KT_n(A) \\
\gamma_n &\colon KU_n(A) \longrightarrow KT_{n-1}(A)  &
\tau_n &\colon KT_n(A) \longrightarrow KO_{n+1}(A)  \; 
\end{align*}
where, for example, the complexification operation $c$ is induced by the inclusion $A \rightarrow \C \otimes A$ and the realification operation $r$ is induced by the inclusion $\C \otimes A \rightarrow M_2(\R) \otimes A$.  For descriptions of the other operations, see Sections~1.1 and 1.2 of \cite{Boer2}.

The target category of united $K$-theory is the category of abstract 
\ct-modules described in \cite{Bou}.  An abstract \ct-module is a triple $M 
= (M\po, M\pu, M\pt)$ consisting of graded modules over $KO_*(\R)$, 
$KU_*(\R)$, and $KT_*(\R)$ respectively.  Furthermore there must be 
$KO_*(\R)$-module homomorphisms $r$, $c$, $\ve$, $\zeta$, $\psi\su$, 
$\psi\st$, $\gamma$ and $\tau$ which satisfy the relations 
\begin{align*}  
rc &= 2    & \psi\su \beta\su &= -\beta\su \psi\su & \xi &= r \beta\su^2 c \\
cr &= 1 + \psi\su & \psi\st \beta\st &= \beta\st \psi\st 
	&  \omega &= \beta\st \gamma \zeta \\ 
r &= \tau \gamma & \varepsilon \beta\so &= \beta\st^2 \varepsilon 
	& \beta\st \varepsilon \tau 
			&= \varepsilon \tau \beta\st + \eta\st \beta\st    \\
c &= \zeta \varepsilon & \zeta \beta\st &= \beta\su^2 \zeta
	 &  \varepsilon r \zeta &= 1 + \psi\st   \\
(\psi\su)^2 &= 1 & \gamma \beta\su^2 &= \beta\st \gamma       
	&  \gamma c \tau &= 1 - \psi\st \\
(\psi\st)^2 &= 1 & \tau \beta\st^2 &= \beta\so \tau & \tau &= -\tau \psi\st  \\
\psi\st \varepsilon &= \varepsilon & 
	\gamma &= \gamma \psi\su 
	\qquad & \tau \beta\st \varepsilon &= 0 \\
\zeta \gamma &= 0 & \eta\so &= \tau \varepsilon 
	& \varepsilon \xi &= 2 \beta\st \varepsilon \\
\zeta &= \psi\su \zeta & \eta\st &= \gamma \beta\su \zeta 
	& \xi \tau &= 2 \tau \beta\st \; 
\end{align*}
as in Section~1.9 of \cite{Bou}.  These relations are satisfied by united $K$-theory by Proposition~1.7 of \cite{Boer2}.  

Not every abstract \ct-module $M$ can be realized as the united $K$-theory 
of a real C*-algebra.  According to Theorem~1.18 of \cite{Boer2} a 
necessary condition is that $M$ be acyclic, i.e. the following complexes 
must be exact:
\begin{align*} 
\cdots \longrightarrow  
&M\pu_{n+1} \xrightarrow{\gamma} 
M\pt_n \xrightarrow{\zeta} 
M\pu_n \xrightarrow{1 - \psi\su}
M\pu_{n} \longrightarrow \cdots  \\
\cdots \longrightarrow  
&M\po_n \xrightarrow{\eta\so} 
M\po_{n+1} \xrightarrow{c} 
M\pu_{n+1} \xrightarrow{r \beta\su^{-1}}
M\po_{n-1} \longrightarrow \cdots  \\
\cdots \longrightarrow  
&M\po_{n} \xrightarrow{\eta\so^2} 
M\po_{n+2} \xrightarrow{\varepsilon} 
M\pt_{n+2} \xrightarrow{\tau \beta\st^{-1}}
M\po_{n-1} \longrightarrow \cdots 
\end{align*}
Our main theorem is that every countable acyclic \ct-module can be realized as the united $K$-theory of a real C*-algebra.  Furthermore, the real C*-algebra can be taken to be simple and purely infinite.
Recall that a complex C*-algebra is said to be a Kirchberg algebra if it is separable, nuclear, simple, and purely infinite.  We say that a real C*-algebra $A$ is a Kirchberg algebra if the complexification $A\sc$ is a Kirchberg algebra.

\begin{thm} ~ \label{main}
\begin{enumerate}
\item[(1)] Let $M$ be any countable acyclic \ct-module.  Then there exists a 
real stable Kirchberg algebra $A$ such that $K\crt(A) \cong M$ and $A\sc$ satisfies the UCT.
\item[(2)] Let $M$ be any countable acyclic \ct-module and let $m$ be any 
element of $M\po_0$ (that is, $m$ is a degree zero element in the real 
part of $M$).  Then there exists a 
real unital Kirchberg algebra $A$ such that $(K\crt(A), [1_A]) \cong (M, m)$ and $A\sc$ satisfies the UCT.
\end{enumerate}
\end{thm}

In \cite{Boer3}, we developed united $KK$-theory  
(generalizing united $K$-theory in the sense that $KK\crt(\R, A) \cong K\crt(A)$ for a real $\sigma$-unital 
C*-algebra $A$) allowing us to state and prove a 
Universal Coefficient Theorem for real C*-algebras.  One of the important 
corollaries highlighting the strength of this theory says that two 
separable C*-algebras $A$ and $B$ such that $A\sc$ and $B\sc$ 
are in the bootstrap category are $KK$-equivalent (in the real sense) 
if and only if $K\crt(A)$ and $K\crt(B)$ are isomorphic \ct-modules.  With 
Theorem~\ref{main} this implies there is an equivalence between the category of countable acyclic \ct-modules and the category of $KK$-equivalence classes of real Kirchberg algebras whose complexification satisfies the UCT.  We believe that this result is a strong indication that united 
$K$-theory should play the same role for real C*-algebras that complex 
$K$-theory plays for complex C*-algebras, especially for any 
classification theorems for real C*-algebras.

Neither real $K$-theory nor complex $K$-theory by itself 
can do the job of united $K$-theory.  In fact, in \cite{Boer2}, we showed 
that the two tensor products of real Cuntz algebras $\mathcal{O}_3\pr 
\otimes \mathcal{O}_3\pr$ and $\mathcal{O}_3\pr \otimes \mathcal{O}_5\pr$ 
are nonisomorphic although their complexifications are isomorphic.  
Many more such examples can be obtained by applying Theorem~\ref{main} (using, for example, the \ct-modules $M_i = \Sigma^i K\crt(\R)$ for $i = 0, 2, 4, 6$).
Hence complex $K$-theory by itself is not sufficient to classify real simple purely infinite
C*-algebras.  Neither is real $K$-theory by itself sufficient, as we will see from 
Theorem~\ref{realnotsuff}.  

In her dissertation \cite{Hewitt}, Beatrice Hewitt  showed that 
acyclic \ct-modules can be classified in terms of their cores, which 
contains only the complex part and the image of $\eta\so$ in the real part 
(and some natural transformations).  Thus the 
self-conjugate part of united $K$-theory is strictly unnecessary.  
However, we know of no way to express \ct~tensor product or $\hom$ functors 
in terms of just the cores.  So for purposes of the K\"unneth formula and 
the universal coefficient theorem, it is still necessary to work with the 
full united $K$-theory.

In our main theorem, we specified the properties that $A\sc$ must satisfy, rather than $A$.  This is partly because it is easier to verify properties in the more familiar setting of complex C*-algebras.  Nothing is lost by this approach since it is usually the case that when $A\sc$ satisfies a certain property, the corresponding property is satisfied by $A$.

For example, any real C*-algebra is simple if its complexification is simple.  
Indeed, if $I$ is a closed ideal in $A$, then $I\sc$ is a closed 
ideal in $A\sc$.  The converse is not true in general.  For 
example, the algebra $\C$ is a simple C*-algebra whose complexification 
$\C \otimes \C \cong \C \oplus \C$ is not simple.  In fact, the 
complexification $A\sc$ is simple if and only if $A$ is simple and 
$A$ is not itself isomorphic to the complexification of a real C*-algebra.

Following the definition in \cite{Stacey}, a real C*-algebra $A$ is purely 
infinite if each hereditary subalgebra of the form $\overline{xAx}$ for a 
nonzero positive element $x$ contains an infinite projection.  Theorem~3.3 
of \cite{Stacey} states that $A$ is purely infinite if $A\sc$ is 
purely infinite.  The converse is still an open question, although there 
is a partial result in Section 4 of \cite{Stacey}.  

We know of no investigations into an intrinsic notion of nuclearity for 
real C*-algebras.  For our purposes in this paper, we define a real 
C*-algebra to be nuclear if and only if its complexification is nuclear.  
At any rate, this is enough to imply that tensor products of real C*-algebras are unique if one of the factors is nuclear.

A real C*-algebra $A$ is said to satisfy the UCT if there is an exact 
sequence of \ct-modules
$$0 \rightarrow {\rm Ext}\scrt(K\crt(A), K\crt(B)) \xrightarrow{\kappa}
	KK\crt(A,B) \xrightarrow{\gamma}
	{\rm Hom}\scrt(K\crt(A), K\crt(B)) 
	\rightarrow 0$$
for all real separable C*-algebras $B$.  According to Theorem~1.1 of 
\cite{Boer3}, a real C*-algebra $A$ satisfies the UCT for united 
$K$-theory if $A\sc$ is in the bootstrap 
category $\mathcal{N}$.  More generally, the proof in that paper shows 
that if $A$ is any real C*-algebra such that $A\sc$ satisfies the 
complex UCT for any complex separable C*-algebra $B$, then $A$ satisfies 
the UCT.  

A partial converse is true since the complex part of the UCT exact sequence for united $K$-theory above is the same as the complex UCT exact sequence for the complex C*-algebras $A\sc$ and $B\sc$.  Thus if a real C*-algebra $A$ satisfies the UCT for all real separable C*-algebras $B$, then the complexification $A\sc$ will satisfy the UCT for all complex separable C*-algebras of the form $B\sc$ for some real C*-algebra $B$.  However, as shown in \cite{Phillips2}, not every complex C*-algebra is the complexification of a real C*-algebra.  

The proof of the main theorem takes place through a series of 
approximating steps.  In Section~\ref{firstconstruction}, we first show 
how to obtain a real separable C*-algebra whose united $K$-theory is isomorphic to the 
prescribed \ct-module.  In Section~\ref{unit} we show how to modify this algebra to form a 
real unital C*-algebra with the same $K$-theory.  Finally in Section~\ref{simple} 
we use a real version of Kumjian's construction to to make our algebra simple and purely infinite, 
allowing us to complete the proof of the main theorem.  

In this paper we will make frequent use of the following important 
theorem, which is implicit in \cite{Boer2} and we state here for 
convenience.  It is an immediate consequence of the results of Section~2.3 
of \cite{Bou} (restated as Propositions~1.14 and 1.15 of \cite{Boer2}) and 
Theorem~1.12 in \cite{Boer2}.

\begin{thm} ~ \label{allforone}
\begin{enumerate} 
\item[(1)] Let $A$ be a real C*-algebra.  If one of the three graded modules 
$KO_*(A)$, $KU_*(A)$, and $KT_*(A)$ is trivial, then all three are trivial.
\item[(2)] Let $f \colon A \rightarrow B$ be a homomorphism of real 
C*-algebras.  If one of the three graded homomorphisms $f_* \colon KO_*(A) 
\rightarrow KO_*(B)$, $f_* \colon KU_*(A) \rightarrow KU_*(B)$, and $f_* 
\colon KT_*(A) \rightarrow KT_*(B)$ is an isomorphism, then all three are 
isomorphisms.
\end{enumerate}
\end{thm}


\section{The First Construction} \label{firstconstruction}

For any acyclic \ct-module $M$, there is according to Theorem~2.9 of \cite{Bou2} a topological spectrum $E$ such that $K\crt(E) \cong M$.  It isn't known in general whether $E$ can be taken to be a actual topological space; however, by Theorem~11.1 of \cite{Bou}, it is possible to find a CW-complex $X$ such that $K\crt(X) \cong M$ if $M$ is finitely generated.  In this section, we prove the following theorem which only requires that $M$ be countable, but leaves the commutative setting far behind.

\begin{thm} \label{surj1} Let $M$ be a countable acyclic \ct-module.  Then 
there is a real separable nuclear C*-algebra $A$ satisfying the UCT such 
that $K\crt(A) \cong M$.
\end{thm} 

First we establish some preliminary notation.  Given a real C*-algebra $A$ 
we define the suspension by
$SA = C_0(\R, A)$ and the desuspension by
$$S^{-1}A  = \{f \in C_0(\R, \C \otimes A) \mid f(-x) =\overline{f(x)} \} 
\; . $$
This nomenclature is justified by the result that $SS^{-1} \R$ (and 
$S^{-1}S\R$)  is $KK$-equivalent to $\R$ (Proposition~1.20 of 
\cite{Boer2}).  More generally, we define
$$S^n A = \begin{cases}
	\underbrace{S S \dots S}_n A & \rm{if}~n \geq 0 \\
	\underbrace{S^{-1} S^{-1} \dots S^{-1}}_{-n} A & \rm{if}~n < 0
		\end{cases}  \; .$$
Let $\iota$ represent the orientation-reversing involution of $SA$, which 
induces multiplication by $-1$ on $K$-theory.		
		
Recall from Section~2.1 of \cite{Boer2} that $K\crt(\R)$, $K\crt(\C)$, and 
$K\crt(T)$ are free \ct-modules.  The \ct-module 
$K\crt(\R)$ is generated by $1\so$, the class of the identity in 
$KO_0(\R)$; the element $\kappa_1 \in KU_0(\C)$ 
generates $K\crt(\C)$ as a \ct-module and satisfies $r(\kappa_1) = 1\su 
\in KO_0(\C)$; and the element $\chi \in KT_{-1}(T)$ 
generates $K\crt(T)$ and satisfies $\tau(\chi) = 1\st \in KO_0(T)$.

\begin{lemma} \label{KOclassspace}
Let $B$ be any real unital C*-algebra and let $x \in KO_0(B)$.  Then there 
is a positive integer $n$ and a C*-algebra homomorphism $\alpha \colon S\R 
\rightarrow M_n S B$ 
such that $\alpha_*(1\so) = x$.
\end{lemma}

Note that we are making use of the identifications $KO_0(\R) = 
KO_{-1}(S\R)$ 
and $KO_0(B) = KO_{-1}(SB)$, claiming that $\alpha_* \colon KO_{-1}(S\R) 
\rightarrow  KO_{-1}(SB)$ sends $1\so$ to $x$.

\begin{proof}
Let $x = [p_1] - [p_2]$ where $p_i$ is a projection in $M_{n_i}(B)$ for $i 
= 1,2$.  First define $\alpha_i \colon \R \rightarrow M_{n_i} B$ by 
$\alpha_i(t) = t p_i$.  Then let $n = n_1 + n_2$ and define $\alpha = S 
\alpha_1 \oplus (S \alpha_2 \circ \iota)$.  Then $\alpha_*(1\so) = 
[p_1] + \iota_*[p_2] = [p_1] - [p_2] = x$.
\end{proof}

\begin{lemma} \label{KUclassspace}
Let $B$ be any real unital C*-algebra and let $y \in KU_0(B)$.  Then there 
is a positive integer $n$ and a
C*-algebra homomorphism $\alpha \colon S\C \rightarrow M_n S^{-1} B$ such 
that 
$\alpha_*(\kappa_1) = \beta\su^{-1}y$.
\end{lemma}

\begin{proof}
Consider the unital inclusion $c \colon \R \rightarrow \C$.  We apply the 
mapping cone construction as in the proof of Theorem~1.18 of \cite{Boer2} 
to obtain a C*-algebra 
homomorphism $\nu \colon S\C \rightarrow 
Cc$.  In that proof, we found that the 
mapping cone $Cc$ is homotopy equivalent to $S^{-1}$ and we proved that 
the element of $KK_{-2}(\C, \R)$ represented by $\nu$ is $\pm r 
\beta\su^{-1}$.  If the sign is negative, replace $\nu$ by $\nu 
\circ \iota$ to make it positive.

Let $y = [p_1] - [p_2]$ where each $p_i$ is a projection in $M_{n_i} \C 
\otimes B$.  
Define a C*-algebra 
homomorphism $\rho_i \colon \C \rightarrow M_{n_i} \C \otimes B$ by 
$\rho_i(t) = t p_i$ for all $t \in \C$.  The 
composition $h_i = \nu \circ S \rho_i$ defines a homomorphism 
from $S\C$ to $M_{n_i} S^{-1} B$ such that $(h_i)_*(1\su) = r 
\beta\su^{-1}[p_i]$.  
Let $n = n_1 + n_2$ and define $h = h_1 \oplus (h_2 
\circ \iota)$ from $S\C$ to $M_n S^{-1} B$, so that $h_*(1\su) 
= r \beta\su^{-1}(y)$.  Then $r h_*(\kappa_1) = h_* r(\kappa_1) = 
h_*(1\su) = r \beta\su^{-1}(y)$.  Since $\ker r = \im 
\beta\su^{-1} c$ (Thoerem~1.18 of \cite{Boer2}), there is an element $x 
\in 
KO_{1}(S^{-1} B) = KO_0(B)$ 
such that $h_*(\kappa_1) = \beta\su^{-1}(y) + \beta\su^{-1} c (x)$.  

To correct the error, let $x = [q_1] - [q_2]$ where $q_i$ is a projection 
in $M_{m_i}B$ for $i = 1,2$.  Define $\mu_i \colon \R \rightarrow 
M_{m_i}B$ 
by $\mu_i(t) = t q_i$ and then define $j_i = S^{-1}\mu_i \circ 
\nu$.  Let $m = m_1 + m_2$ and define $j\colon S \C \rightarrow 
M_mS^{-1}B$ 
by $j = (j_1 \circ \iota) \oplus j_2$.  Since $\beta\su \nu_*(\kappa_1) = 
\nu_* (\beta\su \kappa_1) = r \beta\su^{-1} \beta\su \kappa_1 = r \kappa_1 
= 1\so$, we have 
$\nu_*(\kappa_1) = \beta\su^{-1} c(1\so) \in KU_{-2}(\R)$.  Thus 
$$j_*(\kappa_1) = ((\mu_2)_* - (\mu_1)_*)\beta\su^{-1}c(1\so) = 
	- \beta\su^{-1}c ((\mu_1)_* - (\mu_2)_*)(1\so) = -\beta\su^{-1}c(x) \; .$$

We patch together these two homomorphisms by letting $l = m+n$ and 
defining $\alpha = h \oplus j$ from $S\C$ to $M_l S^{-1}B$ so that 
$\alpha(\kappa_1) = \beta\su^{-1}$y.
\end{proof}

\begin{lemma} \label{KTclassspace}
Let $B$ be any real unital C*-algebra and let $z \in KT_0(B)$.  Then there 
is a 
positive integer $n$ and a C*-algebra homomorphism $\alpha \colon S T 
\rightarrow M_n S^{-2} B$ such that $\alpha_*(\chi) = \beta\st^{-1} z$.
\end{lemma}

\begin{proof}
The mapping cone of the unital inclusion $\ve \colon \R \rightarrow T$ is 
homotopy equivalent to $S^{-2}$ (as in the proof of Theorem~1.18 of 
\cite{Boer2}).  Thus
 we obtain a C*-algebra homomorphism $\sigma \colon ST \rightarrow 
S^{-2}$.  Also in the proof of Theorem~1.18, we proved that 
the element of $KK_{-3}(T, \R)$ represented by $\sigma$ is $\pm \tau 
\beta\st^{-1}$.  If the sign is negative, replace $\sigma$ by $\sigma 
\circ \iota$ to make it positive.

Let $z = [p_1] - [p_2]$ where $p_i$ is a projection in $M_{n_i} T \otimes 
B$ for $i = 1,2$.  Since $T$ is commutative, there is a C*-algebra 
homomorphism $\rho_i \colon T \rightarrow M_{n_i} T \otimes B$ defined by 
$\rho_i(t) = t p_i$ for all $t \in T$.  The 
composition $h_i = \sigma \circ S \rho_i$ defines a homomorphism 
from $ST$ to $M_{n_i} S^{-2} B$ such that $(h_i)_*(1\st) = \tau 
\beta\st^{-1}[p_i]$.  
Let $n = n_1 + n_2$ and define $h = h_1 \oplus (h_2 \circ \iota)$ from 
$ST$ to $M_n S^{-2} B$ so that $h_*(1\st) 
= \tau \beta\st^{-1}(z)$.  Then $\tau h_*(\chi) = h_* \tau(\chi) = 
h_*(1\st) = \tau \beta\st^{-1}(z)$.  Since $\ker \tau = \im 
\beta\st^{-1} \ve$ (Theorem~1.18 of \cite{Boer2}), there is an element $x 
\in 
KO_{2}(S^{-2} B) = KO_0(B)$ 
such that $h_*(\chi) = \beta\st^{-1}(z) + \beta\st^{-1} \ve(x)$.  

To correct the error, let $x = [q_1] - [q_2]$ where $q_i$ is a projection 
in $M_{m_i}B$ for $i = 1,2$.  Define $\mu_i \colon \R \rightarrow 
M_{m_i} B$ by $\mu_i(t) = t q_i$ and then define $j_i = S^{-2} 
\mu_i \circ \sigma$.  Let $m = m_1 + m_2$ and define $j \colon ST 
\rightarrow M_mS^{-2} B$ by $j = (j_1 \circ \iota) \oplus j_2$.  Since 
$\beta\st \sigma_*(\chi) = \sigma_*(\beta\st(\chi) = \tau (\chi) = 1\st$, 
we have $\sigma_*(\chi) = \beta\st^{-1} (1\st)$ in $KT_{-4}(\R)$.  Thus
$$j_*(\chi) 
	= ((\mu_2)_* - (\mu_1)_*) \beta\st^{-1}\ve(1\so)
	= - \beta\st^{-1} \ve ((\mu_1)_* - (\mu_2)_*) (1\so) 
	= - \beta\st^{-1} \ve(x) \; .$$

We patch these two homomorphisms together by letting $l = m+n$ and 
defining $\alpha = h \oplus j$ from $ST$ to $M_l S^{-2 }B$ so that 
$\alpha(\chi) = \beta\st^{-1}(z)$.
\end{proof}

\begin{proof} [Proof of Theorem~\ref{surj1}]
If $M$ is a free \ct-module, then it can be written as a direct sum of 
monogenic free \ct-modules, and each monogenic \ct-module can be realized 
as the united $K$-theory of $\R$, $\C$, $T$, or a suspension thereof.  
Therefore, $M$ can be realized as the united $K$-theory of a direct sum of 
countably many such C*-algebras.

Now, let $M$ be an arbitrary countable acyclic \ct-module.  By 
Theorems~3.2 
and 3.4 in \cite{Bou}, we can find a resolution
$$ 0 \rightarrow F_1 \xrightarrow{\mu_1} F_0 \xrightarrow{\mu_0} M 
			\rightarrow 0 \; $$
where $F_0$ and $F_1$ are countable and free \ct-modules.

As in the first paragraph, find real separable C*-algebras $B$ and $C$ 
such that 
$F_0 = K\crt(B)$ 
and $F_1 = K\crt(C)$.  In particular, we set 
$$B = \bigoplus_{i \in I\so} S^{k_i} \R \oplus \bigoplus_{i \in 
I\su} S^{k_i} \C \oplus \bigoplus_{i \in I\st} S^{k_i} T$$ 
where $I\so$, $I\su$, and $I\st$ are disjoint countable index sets and 
where  
$k_i \in \{0, 1, \dots, 7 \}$ for each $i$.  Our strategy is to realize 
$\mu_1$ geometrically.  That is, we wish to produce a C*-algebra 
homomorphism $\beta \colon B \rightarrow C$ whose induced homomorphism on 
united $K$-theory is $\mu_1$.  Actually, we will replace $B$ and $C$ with 
algebras $B'$ and $C'$ and the induced homomorphism $\beta_* \colon 
K\crt(B') \rightarrow K\crt(C')$ will not be identical to $\mu_1$ but will 
be injective and will have the same cokernel as $\mu_1$.

For any unital C*-algebra $D$, let $S^{\sim-} D = (S^{-1} D)^{\sim}$ 
denote 
the unitized desuspension of $D$ and let $S^{\sim -n} D$ denote the 
$n$-fold unitized desuspension.  Let $C'' = S^{10} S^{-2} S^{\sim -8} 
(C^\sim)$ and let $C' = \mathcal{K} \otimes C''$ where $\mathcal{K}$ is an 
algebra of compact operators on a separable Hilbert space.
Note that $K\crt(C') = K\crt(C) \oplus K\crt(S^{10} S^{-2} S^{\sim - 8} 
\R)$ because of the split exact sequence
$$0 \rightarrow S^{10} S^{-2} S^{\sim - 8} \R
	\rightarrow S^{10} S^{-2} S^{\sim - 8} C^{\sim}
	\rightarrow S^{10} S^{-10} C
	\rightarrow 0 \; .$$

For each $i \in I\so$, we construct a geometric realization of the 
restricted homomorphism
$$K\crt(S^{k_i} \R) \rightarrow K\crt(C) \; $$
as follows.  Let $x 
\in KO_{-k_i}(C) = KO_0(S^{-k_i} C)$ be the image of $1\so 
\in KO_{-k_i}(S^{k_i} \R) = KO_0(\R)$.  By
Lemma~\ref{KOclassspace} there is a homomorphism
$\alpha_i \colon S \R \rightarrow M_{n_i} S (S^{-k_i}  C)^{\sim}$ such 
that 
$(\alpha_i)_*(1\so) = x$.  
Then apply the suspension and desuspension operations to $\alpha_i$ and 
follow it by the inclusion into $C''$ to form the homomorphism 
$$\beta_i \colon S^{10} S^{k_i - 10} \R 
	\rightarrow M_{n_i} S^{10} S^{k_i-10} (S^{-k_i} C)^{\sim} 
	\hookrightarrow M_{n_i} C'' \; $$
which agrees on united $K$-theory with the restriction of 
$\mu_1$ to $K\crt(S^{k_i} \R)$.

Similarly, for each $i \in I\su$, consdider the restriction of $\mu_1$
$$K\crt(S^{k_i}\C) \rightarrow K\crt(C) \; ,$$
and let $y_i \in KU_{-k}(C)$ be the image of $\kappa_1 \in 
KU_{-k_i}(S^{k_i} \C)$.  Using Lemma~\ref{KUclassspace}, let
$\alpha_i \colon S \C \rightarrow M_{n_i} S^{-1} (S^{-k_i} C)^{\sim}$ be 
given 
satisfying 
$(\alpha_i)_*(\kappa_1) = \beta\su^{-1} y$.  
Again suspend and desuspend to form the composition
$$\beta_i \colon S^{11} S^{k_i - 9} \C 
	\rightarrow M_{n_i} S^{10} S^{k_i-10} (S^{-k_i} C)^{\sim}
	\hookrightarrow M_{n_i} C'' \; .$$
The induced homomorphism $(\beta_i)_*$ on united $K$-theory agrees with 
the restriction of 
$\mu_1$ to $K\crt(S^{k_i} \C)$ up to multiplication by $\beta\su^{-1}$.  
This is not a problem for us;  since $\beta\su^{-1}$ is an isomorphism on 
united $K$-theory, the homomorphism $(\beta_i)_*$ is still injective and 
its image is the same as that of $\mu_1$.

Thirdly, for each $i \in I\su$, consider the restriction of $\mu_1$
$$K\crt(S^{k_i} T) \rightarrow K\crt(C)$$
and let $z_i \in KT_{-k_i-1}(C)$ be the image of $\chi \in 
KT_{-k_i-1}(S^{k_i} \C)$.  By Lemma~\ref{KTclassspace}, let
$\alpha_i \colon ST \rightarrow M_{n_i} S^{-2} (S^{-k_i-1}  C)^{\sim}$ 
be given satisfying 
$(\alpha_i)_*(\chi) = \beta\st^{-1} z$.  
Again suspend and desuspend to form 
$$\beta_i \colon S^{11} S^{k_i - 7} T 
	\rightarrow M_{n_i} S^{10} S^{k_i - 9}  (S^{-k_i - 1} C)^{\sim}
	\hookrightarrow M_{n_i} C''	 \; , $$
a map which on united $K$-theory agrees with the restriction of 
$\mu_1$ to $K\crt(S^{k_i} \C)$ up to multiplication by $\beta\st^{-1}$.

We need one more homomorphism, 
$$\beta_0 \colon S^{10} S^{-2} S^{\sim -8} \R
	\rightarrow S^{10} S^{-2} S^{\sim -8} (C^{\sim})$$
based on the unital inclusion $\R \rightarrow C^{\sim}$.

We assemble the homomorphisms using a big Hilbert space.  Let 
$\mathcal{K}$ be the algebra of compact operators on
a separable Hilbert space and let $\phi_i$ be a collection of mutually 
orthogonal 
inclusions from $M_{n_i}$ to $\mathcal{K}$ for $i \in I\so \cup I\su 
\cup I\st \cup \{0\}$.  Let 
$$B' = S^{10} S^{-2} S^{\sim -8} \R \oplus
		\bigoplus_{i \in I\so} S^{10} S^{k_i-10} \R \oplus 
		\bigoplus_{i \in I\su} S^{11} S^{k_i-9} \C \oplus	
		\bigoplus_{i \in I\st} S^{11} S^{k_i-7} T  
		 $$  
and we define $\beta \colon B' \rightarrow C' = \mathcal{K} \otimes C''$ 
by setting it to be $\phi_i \circ \beta_i$ on each summand.  

Therefore, we have 
a geometric realization of $\mu_1$ in the sense that $\beta_*$ is 
injective and has 
the same cokernel as $\mu_1$.  Let $A'$ be the mapping cone of $\beta$.  
Then we have a short exact sequence
$$0 \rightarrow SC' \rightarrow A' \rightarrow B' \rightarrow 0 \; .$$
In the resulting long exact sequence, the homomorphism $K\crt(B') 
\rightarrow K\crt(SC')$ of degree $-1$ is the same as $\beta_*$ (see 
Proposition~2.5 of \cite{Schochet3} or Theorem~1.1 of \cite{CS}).  
Since $\beta_*$ is injective, 
the long exact sequence collapses to the short exact sequence
$$0 \rightarrow K\crt(B') \xrightarrow{\beta_*}
		K\crt(C') \xrightarrow{i_*}
		K\crt(A') \rightarrow 0$$		
where $i_*$ has degree $-1$.  The united $K$-theory of $A'$ 
is thus a shift of the \ct-module $M$ so the algebra $A = S^{-1}A'$ 
finishes the 
job.

Since $A$ is constructed from the commutative algebras $\R$, $\C$ and $T$ 
using the 
operations of countable direct sum, suspensions, desuspensions, 
unitization, forming matrix algebras, stabilization, and forming mapping 
cones
we know that $A$ is separable, nuclear, 
and in the category of real C*-algebras that satisfy the Universal 
Coefficient 
Theorem.
\end{proof}


\section{Unital} \label{unit}

The goal of this section is to show that given a real C*-algebra $A$, we 
can obtain a unital algebra with the same united 
$K$-theory.  For this, we will use the real analog of the construction of 
Proposition~4.1 in \cite{BKP}.  

We begin by recording some results regarding real simple purely infinite C*-algebras and their $K$-theory.  These results are analogs of well-known results in the theory of complex simple purely infinite C*-algebras.  In each case, the proof follows directly from the corresponding result in the complex case, or can be proven in the same way as the complex version.

It is well-known that the inclusion of a full corner in a complex 
C*-algebra induces an isomorphism on $K$-theory.  It is an easy 
consequence of Theorem~\ref{allforone} that the same is true for real 
C*-algebras.  For completeness, we record the proofs of both statements 
below.

\begin{prop} ~ \label{fullcorner}
\begin{enumerate} 
\item[(1)] Let $p$ be a full projection in a complex C*-algebra $A$.  
Then the inclusion $i \colon pAp \rightarrow A$ induces an isomorphism on 
$K$-theory.
\item[(2)] Let $p$ be a full projection in a real C*-algebra $A$.  Then 
the inclusion $i \colon pAp \rightarrow A$ induces an isomorphism on 
united $K$-theory.
\end{enumerate} \end{prop}

\begin{proof}
Let $A$ be a complex C*-algebra and let $p$ be a full projection.  By 
Lemma~2.5 of \cite{Brown}, there is a partial isometry $v \in M(pAp 
\otimes \mathcal{K})$ such that $v^* v = 1$ and $v v^* = p \otimes 1$.  
Replacing $v$ by $(p \otimes 1)v$, we may assume that $v \in (pAp \otimes 
\mathcal{K})^+$.  Then there is an isomorphism $\alpha \colon pAp \otimes 
\mathcal{K} 
\rightarrow A \otimes \mathcal{K}$ defined by $x \mapsto v^* x v$.  

Now, if $q$ is any projection in $(pAp \otimes \mathcal{K})^+$, then 
$(qv)^*(qv) = v^* q v$ and $(qv)(qv)^* = q (p \otimes 1) q^* = q$.
Thus in $K_0(A^+)$ we have $[i(q)] = [q] = [v^* q v] = [\alpha(q)]$.
Similarly, if $u$ is any unitary in $(pAp \otimes \mathcal{K})^+$, then 
in $K_1(A^+)$ we have $[i(u)] = [u] = [v^* u v] = [\alpha(u)]$.

Therefore, $i_*$ and $\alpha_*$ agree as homomorphisms from $K_*(pAp)$ to 
$K_*(A)$.  Since $\alpha_*$ is an isomorphism, so is $i_*$.  This proves 
part 
(1).
To prove part (2), let $p$ be a full projection in a real C*-algebra $A$.  
By part (1) the inclusion $i_*$ induces an isomorphism on complex 
$K$-theory $KU_*(pAp) \rightarrow KU_*(A)$.  Therefore, $i_*$ is an 
isomorphism 
on united $K$-theory by Theorem~\ref{allforone}.
\end{proof}

\begin{lemma} \label{realpi}
Let $p$ and $q$ be non-trivial projections in a simple purely infinite 
C*-algebra.  Then there is a projection $p'$ such that $p' \sim p$ and $p' 
< q$.
\end{lemma}

The complex version of Lemma~\ref{realpi} can be 
found as Proposition~1.5 in \cite{Cuntz} or Lemma~V.5.4 in 
\cite{Davidson}.  The proof of Lemma~\ref{realpi} follows exactly the proof 
of Lemma~V.5.4 in \cite{Davidson}.  (This was also observed by Stacey in the proof of Proposition~4.1 in \cite{Stacey}.  Once this lemma is established, the 
proof of Proposition~\ref{realc} below follows exactly the proof of Theorem~1.4 
of \cite{Cuntz}.

\begin{prop} \label{realc}
Let $A$ be a real simple purely infinite C*-algebra.  Then
$$KO_0(A) \cong \{ [p]  \mid p \text{ is a non-zero projection in } A \}$$
where $[p]$ represents the Murray-von Neumann equivalence class of a 
projection $p$ in $A$.
\end{prop}

\begin{prop} \label{unitfunctor}
There is a functor $F$ from the category of all real C*-algebras (and 
real C*-algebra homomorphisms) to the 
category of all real unital C*-algebras (and real unital C*-algebra 
homomorphisms) and a natural transformation $\eta \colon A 
\rightarrow F(A)$ which induces an isomorphism on united $K$-theory.  
Furthermore,
\begin{enumerate}
\item[(1)] If $A$ is nuclear, then $F(A)$ is nuclear.
\item[(2)] If $A$ is separable, then $F(A)$ is separable and $\eta$ is a 
$KK$-equivalence.
\item[(3)] If $A$ is separable and satisfies the UCT, then $F(A)$ satisfies the 
UCT.
\end{enumerate}
\end{prop}

\begin{proof}
Let $\mathcal{O}_{\infty}\pr$ be the real 
Cuntz algebra generated by a sequence of mutually orthogonal isometries.  
By Theorem~\ref{allforone} the unital inclusion $\R \rightarrow \mathcal{O}_{\infty}\pr$ induces an 
isomorphism on united $K$-theory since the complexification
$\C \rightarrow \mathcal{O}_{\infty}$ induces an isomorphism on $K$-theory.
By Proposition~\ref{realc}, there is a non-zero 
projection $e \in \mathcal{O}_{\infty}\pr$ and a projection $q < e$ such 
that $[e] = 0$ and $[q] = [1_{\mathcal{O}_\infty\pr}]$.

Since $e$ is infinite, there exists a proper subprojection $p_1$ such that 
$p_1 \sim e$.  Let $p_2 = e - p_1$.  Then $[p_1] = [p_2] = [e] = 0$.  Therefore 
(again by Proposition~\ref{realc}) there are partial isometries 
$s_1$ and $s_2$ in $e \mathcal{O}_\infty\pr e$
such that $s^*_i s_i = e$ and $s_i s^*_i = p_i$.  Let $D = C^*(s_1, s_2)$.  
Then the algebra $D = C^*(s_1, s_2)$ is a unital subalgebra of 
$e \mathcal{O}_{\infty}\pr e$ which is isomorphic to $\mathcal{O}_2\pr$.

Now, for any real C*-algebra $A$, let 
$A^+$ be the unitization of $A$ and let $\pi_A \colon A^+ \rightarrow 
\R$ be the usual projection with kernel $A$.  We define 
$$F(A) = \{b \in e \mathcal{O}_{\infty}\pr e \otimes A^+ \mid (1 \otimes 
\pi_A)(b) \in D \} \; .$$
The element $e \otimes 1$ is a unit for $F(A)$.  The natural 
transformation $\eta \colon A \rightarrow F(A)$ is defined by $a \mapsto 
q \otimes a$.

We will show 
that $\eta$ induces an isomorphism on united $K$-theory.  Note that 
$\eta$ is a composition of the homomorphism $A \rightarrow e 
\mathcal{O}\pr_\infty e \otimes A$ defined by $a \mapsto q \otimes a$ and 
the 
inclusion $e \mathcal{O}\pr_\infty e \otimes A \hookrightarrow F(A)$.  
The homomorphism 
$A \rightarrow e \mathcal{O}\pr_\infty e \otimes A$ 
induces an isomorphism on united $K$-theory because 
the map $\R \rightarrow e \mathcal{O}\pr_\infty e$ defined by $t \mapsto 
tq$ does using the K\"unneth formula for united $K$-theory (\cite{Boer2}).  
Secondly, the inclusion $e \mathcal{O}\pr_\infty e \otimes A 
\hookrightarrow F(A)$ induces an isomorphism on united $K$-theory because 
of the short exact sequence
$$
0 \rightarrow e \mathcal{O}\pr_{\infty}e \otimes A
	\hookrightarrow F(A) \xrightarrow{1 \otimes \pi_A} D \rightarrow 0
$$
and the fact that $K\crt(D) = K\crt(\mathcal{O}_2) = 0$.  
It follows that $\eta$ induces 
an isomorphism on united $K$-theory.

It is clear from the short exact sequence above that if $A$ is separable or 
nuclear, then the same is true of $F(A)$.  Furthermore, the argument of 
the previous paragraph also works for $KK$-theory, showing that $\eta$ 
induces isomorphisms
$$KK\crt(B, A) \rightarrow KK\crt(B, F(A))$$
and
$$KK\crt(F(A), B) \rightarrow KK\crt(A, B)$$
for any real separable C*-algebra $B$.  If $A$ 
is separable, then so is $F(A)$ and by the Yoneda Lemma, 
$\eta\eta$ induces a $KK$-equivalence.  In particular, 
if $A$ is separable and satisfies the UCT, so does $F(A)$.
\end{proof}


\section{Simple and Purely Infinite} \label{simple}

In \cite{Kumjian} Alex Kumjian presents a construction (based on a 
special case of Michael Pimsner's construction in \cite{Pimsner}) which turns any complex 
separable unital C*-algebra $A$ into a complex C*-algebra $\mathcal{O}_E$ 
which is simple and purely infinite such that there is an inclusion $A 
\hookrightarrow \mathcal{O}_E$ which is a (complex) $KK$-equivalence.  In 
this section, we show that this construction can be applied to the real 
case.  Combined with the results from Sections~\ref{firstconstruction} and 
\ref{unit}, this will complete the proof of Theorem~\ref{main}.

\begin{prop} \label{realkumjian}
Let $A$ be a real separable unital C*-algebra.  Then there is a real 
separable simple purely infinite C*-algebra $\mathcal{O}\pr_E$ and a 
unital inclusion $A \hookrightarrow \mathcal{O}\pr_E$ which induces an 
isomorphism on united $K$-theory.  Furthermore, if $A$ is nuclear and 
satisfies the UCT, then the same is true of $\mathcal{O}\pr_E$ and $\iota$ 
is a (real) $KK$-equivalence.
\end{prop}

Recall that a complex C*-algebra $A$ is said to have a real structure if 
there is a conjugate linear involution $x \mapsto \overline{x}$.  In that 
case, the set $A\pr$ of fixed points is a real C*-algebra.  Conversely, 
given a real C*-algebra $A$, the complexification $A\sc$ 
has a real structure given by $a_1 + i a_2 \mapsto a_1 - i a_2$.  These 
functors are inverse to each other so there is a bijection between complex 
C*-algebras with real structure and real C*-algebras.  To prove 
Proposition~ \ref{realkumjian} we will retrace Kumjian's construction, 
showing that the real structure of $A\sc$ passes to $\mathcal{O}_E$.

\begin{defn}
Let $A$ be a complex C*-algebra with a real structure.  
\begin{enumerate}
\item[(1)]
A Hilbert $A$-module $E$ is said to have a real structure if $E$ has a 
conjugate linear involution $e \mapsto \overline{e}$ that satisfies
$\overline{ \langle e, f \rangle} = \langle \overline{e}, \overline{f} 
\rangle$ and
$\overline{e \cdot a} = \overline{e} \cdot \overline{a}$ for all $a \in A$ 
and $e,f \in E$.
\item[(2)] A Hilbert $A$-bimodule $(E, \phi)$ is said to have a real structure 
if the Hilbert $A$-module $E$ has a real structure as in part (1) and the 
homomorphism $\phi \colon A \rightarrow \mathcal{L}(E)$ satisfies 
$\overline{\phi(a) e} = \phi(\overline{a}) \overline{e}$ for all $a \in A$ 
and $e \in E$.
\end{enumerate}
\end{defn}

If $E$ is a Hilbert $A$-module with a real structure, then the C*-algebra 
$\mathcal{L}(E)$ has a real structure defined by $\overline{T}(e) = 
\overline{T(\overline{e})}$ for all $T \in \mathcal{L}(E)$ and $e \in E$.  
With this language, the Hilbert bimodule condition above can be restated as
$\overline{\phi(a)}  = \phi (\overline{a}) $, saying that the 
*-homomorphism $\phi \colon A \rightarrow \mathcal{L}(E)$ respects the 
real structures.

Kumjian's construction begins with a faithful representation $\pi \colon A 
\rightarrow \mathcal{L}(H)$ where $H$ is a separable complex Hilbert space 
such that $\pi(A) \cap \mathcal{K}(H) = \{ 0 \}$.  If $A$ has a real 
structure, we can start with a representation of $A\sr$ on a real Hilbert 
space $H\sr$ and then complexify.  Thus we can assume that $\pi$ respects 
the real structures of $A$ and $\mathcal{L}(H)$.  Following Kumjian, we 
define a Hilbert $A$-bimodule $(E, \phi)$ by 
$$E = H \otimes\sc A $$
with bimodule structure given by
$(\xi \otimes a) \cdot b = \xi \otimes (a \cdot b)$
and $\phi(b)(\xi \otimes a) = \pi(b) \xi \otimes a$ for all $a,b \in A$ 
and $\xi \in H$.  We give $(E, \phi)$ a real structure by $\overline{\xi 
\otimes a} = \overline{\xi} \otimes \overline{a}$.

Similarly, the Fock space 
$$\mathcal{E}_+ = \bigoplus_{n=0}^\infty E^{\otimes n}$$ 
is also a Hilbert $A$-bimodule with a real structure.  The involution is 
defined on pure tensors by
$$\overline{e_1 \otimes e_2 \otimes \dots \otimes e_n} =
 \overline{e_1} \otimes \overline{e_2} \otimes \dots \otimes 
\overline{e_n} \; .$$
For any element $e \in E$, we define the operator $T_e \in 
\mathcal{L}(\mathcal{E}_+)$ on pure tensors by $T_e( e_1 \otimes \dots 
\otimes e_n) = e \otimes  e_1 \otimes \dots \otimes e_n$.  Since 
$\overline{T_e} = T_{\overline{e}}$, the involution of $\mathcal{L}(E)$ 
restricts to an involution of the algebra $\mathcal{T}_E$ generated by 
$\{T_e\}_{e \in E}$.

In the general case, $\mathcal{O}_E$ is the quotient of $\mathcal{T}_E$ by 
the C*-algebra generated in $\mathcal{L}(\mathcal{E}_+)$ by 
$\mathcal{L}\left(\bigoplus_{n=0}^N E^{\otimes n}\right)$ for all positive 
integers $N$.  But under the assumption $\pi(A) \cap \mathcal{K}(H) = \{0 
\}$, we have $\mathcal{O}_E \cong \mathcal{T}_E$ (see \cite{Pimsner}, 
Corollary~3.14).  In either case, the involution of 
$\mathcal{L}(\mathcal{E}_+)$ induces one on $\mathcal{O}_E$.
Furthermore, the inclusion $\iota \colon A \hookrightarrow \mathcal{O}_E$ 
given by $\iota(a)(e_1 \otimes e_2 \otimes \dots \otimes e_n) = 
\phi(a)(e_1) \otimes e_2 \otimes \dots \otimes e_n$, respects the real 
structures of $A$ and $\mathcal{O}_E$.

If we begin with a real separable unital C*-algebra $A$, then the 
complexification $A\sc$ has a real structure and the construction above 
yields an inclusion $\iota \colon A \rightarrow \mathcal{O}_E\pr$ where 
$\mathcal{O}_E\pr$ is the fixed point set of $\mathcal{O}_E$.

\begin{proof}[Proof of Proposition~\ref{realkumjian}]
Let $A$ be a real separable unital C*-algebra.  Applying the construction 
above, we obtain an inclusion $\iota \colon A \rightarrow 
\mathcal{O}_E\pr$.  By Theorem~2.8 of \cite{Kumjian}, $\mathcal{O}_E$ is 
simple and purely infinite.  Thus $\mathcal{O}_E\pr$ is simple and purely 
infinite by Theorem~3.3 of \cite{Stacey}.  By Corollary~4.5 of 
\cite{Pimsner}, the inclusion $\iota \colon A\sc \rightarrow 
\mathcal{O}_E$ is a $KK$-equivalence.  In particular, it induces an 
isomorphism on $K$-theory, so by Theorem~\ref{allforone}, $\iota \colon A 
\rightarrow \mathcal{O}_E\pr$ induces an isomorphism on united $K$-theory.

If $A$ is nuclear and satisfies the UCT, then by Theorem~3.1 of 
\cite{Kumjian}, the same is true of $\mathcal{O}_E$.  Thus 
$\mathcal{O}_E\pr$ is nuclear and satisfies the UCT.  In particular, since 
$A$ and $\mathcal{O}_E\pr$  have isomorphic united $K$-theory and  both 
satisfy the UCT, they are $KK$-equivalent.
\end{proof}

The proof of Section~4 of \cite{Pimsner} will probably carry over to show 
that $\iota$ is a (real) $KK$-equivalence in general, giving a stronger 
statement than our Theorem~\ref{realkumjian}, but we don't need this for 
our present purposes.

Note that absent a full classification theorem for real simple purely infinite 
C*-algebras, there is no guarantee that $\mathcal{O}\pr_E$ is independent 
of the choice of $\pi$ (as $\mathcal{O}_E$ is when $A\sc$ is nuclear and 
satisfies the UCT).

\begin{proof}[Proof of Theorem~\ref{main}]
Let $M$ be a countable acyclic \ct-module.  By Theorem~\ref{surj1}, there
is a real separable nuclear C*-algebra $A_1$ satisfying the UCT such that
$K\crt(A_1) \cong M$.  Applying the functor $F$ of
Proposition~\ref{unitfunctor}, there is a real separable nuclear unital
C*-algebra $A_2$ satisfying the UCT such that $K\crt(A_2) \cong M$.  Then
applying the real Kumjian construction (Proposition~\ref{realkumjian}),
there is a real separable nuclear unital simple purely infinite $A_3$
satisfying the UCT such that $K\crt(A_3) \cong M$.  Finally, let $A_4 =
\mathcal{K}(H) \otimes\sr A_3$ where $H$ is a real separable Hilbert space. 
By Lemma~\ref{realstable} below, $A_4$ is purely infinite and is the real 
C*-algebra needed to prove part (1).

Now, let $m$ be any element in $M_0\po$.  By Proposition~\ref{realc} there 
is a projection $p \in A_4$ such that $[p] = m$.  Let $A_5$ be the corner 
algebra $pA_4p$.  By Proposition~\ref{fullcorner}, $K\crt(A_5) \cong M$.  
This proves part (2). 
\end{proof}

\begin{lemma} \label{realstable}
If $A$ is a real purely infinite simple C*-algebra, then the stabilization 
$\mathcal{K}(H) \otimes\sr A$ is also purely infinite and simple.
\end{lemma}

\begin{proof}
By Lemma~4.2 of \cite{Stacey}, the matrix algebras $M_n(A)$ are purely 
infinite.  The proof of Proposition~4.1.8 of \cite{Rordam} carries over 
immediately to the real case to show that the inductive limit of simple purely infinite C*-algebras is 
again simple and purely infinite.
\end{proof}


\section{The inadequacy of real $K$-theory} \label{example}

In this section, we will draw one small application from our main theorem, creating an example which shows that $K$-theory by itself cannot classify isomorphism classes or even $KK$-equivalence of real simple purely infinite C*-algebras.  Further applications of Theorem~\ref{main} will appear in \cite{Boer5}.

\begin{thm} \label{realnotsuff}
There exist two real C*-algebras $A$ and $B$ such that $KO_*(A) \cong KO_*(B)$, but $K\crt(A) \ncong K\crt(B)$.
\end{thm}

\begin{proof}
By Theorem~\ref{main}, it suffices to find two distinct countable acyclic \ct-modules whose real parts are isomorphic.  I am indebted to A.K. Bousfield for sharing with me the example of such \ct-modules. We will employ a \ct-module construction found in Chapter~8 of \cite{Hewitt}.   

Let $(G, \alpha)$ be a group with involution satisfying $\ker(1+\alpha) = \im(1-\alpha)$ and $\ker(1-\alpha) = \im(1+\alpha)$.  Let $G^+ = \{ g \in G \mid \alpha(g) = g \}$ and $G^- = \{g \in G \mid \alpha(g) = -g \}$.  Then there are exact sequences
$$0 \rightarrow G^+ \xrightarrow{i^+} G \xrightarrow{\pi^-} G^- \rightarrow 0$$
and  
$$0 \rightarrow G^- \xrightarrow{i^-} G \xrightarrow{\pi^+} G^+\rightarrow 0$$
where $\pi^{+} = 1 + \alpha$, $\pi^- = 1 - \alpha$, and $i^+$ and $i^-$ are the inclusion homomorphisms.  Then Table~\ref{P} displays the groups and natural transformations of an acyclic \ct-module $P(G, \alpha)$.  It is easy, if tedious, to verify that the \ct relations hold and that the sequences are exact making it acyclic.

\begin{table}[h]
\caption{$P(G,\alpha)$} \label{P}
$$\begin{array}{|c|c|c|c|c|c|c|c|c|c|}  
\hline \hline  
n & 0 & 1 & 2 & 3 & 4 & 5 & 6 & 7 & 8 \\
\hline  \hline
N(G)\po
& G^+ & 0 & G^- & 0 & G^+ & 0 & G^- & 0  & G^+ \\
\hline  
N(G)\pu 
& G & 0 & G & 0 & G & 0 & G & 0 & G  \\
\hline  
N(G)\pt
& G^+ & G^- & G^- & G^+ & G^+ & G^- & G^- & G^+ & G^+ \\
\hline \hline
(\eta\so)_n & 0&0&0&0&0&0&0&0&0    \\
\hline
c_n & i^+ & 0 & i^- & 0 & i^+ & 0 & i^- & 0 & i^+    \\
\hline
r_n & \pi^+ & 0 & \pi^- & 0  & \pi^+ & 0 & \pi^- & 0 & \pi^+    \\
\hline
\varepsilon_n & 1 & 0 &  1 & 0 & 1 & 0 & 1 & 0 & 1  \\
\hline
\zeta_n &  i^+ & 0 & i^- & 0 &  i^+ & 0 & i^- & 0 & i^+     \\
\hline
(\psi\su)_n & \alpha & 0 & -\alpha & 0 & \alpha & 0 & -\alpha & 0 & \alpha    \\
\hline
(\psi\st)_n & 1 & -1  & 1 & -1 & 1 & -1 & 1 & -1 & 1  \\
\hline
\gamma_n & \pi^+ & 0 & \pi^- & 0 & \pi^+ & 0 & \pi^- & 0 &  \pi^+    \\
\hline
\tau_n & 0 & 1 &  0 & 1 & 0 & 1 & 0 & 1 & 0      \\
\hline \hline
\end{array}$$
\end{table}

Let $G = \Z_2^4$ and $H = \Z_4 \oplus \Z_2^2$ with involutions 
$$\alpha = \left( \begin{matrix}  0 & 1 & 0 & 0 \\ 1 & 0 & 0 & 0 \\ 0 & 0 & 0 & 1 \\ 0 & 0 & 1 & 0 
	\end{matrix} \right) \qquad \text{and} \qquad
    \beta = \left( \begin{matrix} 1 & 0 & 2 \\ 1 & 1 & 0 \\ 0 & 0 & 1
	\end{matrix} \right) \;  $$
respectively.
Then $G^+ \cong H^+ \cong \Z_2^2$ and $G^- \cong H^- \cong \Z_2^2$.  Thus the real parts of $P(G, \alpha)$ and $P(H, \beta)$ agree while the complex parts do not.
\end{proof}

\end{document}